\documentstyle [12pt] {article}
\input{amssym.def}
\input{amssym.tex}

\newcommand{\N}{{\Bbb N}}
\newcommand{\Q}{{\Bbb Q}}
\newcommand{\R}{{\Bbb R}}
\newcommand{\Z}{{\Bbb Z}}

\newcommand{\ga}{{\gamma}}
\newcommand{\de}{{\delta}}
\newcommand{\De}{{\Delta}}
\newcommand{\la}{{\lambda}}
\newcommand{\La}{{\Lambda}}
\newcommand{\Om}{{\Omega}}
\newcommand{\pfend}{$/\!\! /$\bigskip}

\newcommand{\lex}{\times\hspace{-.169in}\raisebox{5pt}{$\scriptsize
\leftarrow$}\, }

\newtheorem{defi}{Definition}[section]%master template; adds %section
\newtheorem{theorem}{Theorem}

\newtheorem{lemma}[defi]{Lemma}
\newtheorem{proposition}[defi]{Proposition}

\newtheorem{corollary}[defi]{Corollary}

\newtheorem{*example}[defi]{*Example}

\addtocontents{ences.txt}{References}

\begin{document}

\title{\bf  Finitely generated lattice-ordered groups with soluble word problem.}

\author{\sc A. M. W. Glass}

\small{\date{Submitted 19$^{th}$ March 2007}}

\maketitle

\begin{abstract}
{William W. Boone and Graham Higman proved that a finitely generated group has soluble word problem if and only if it can be embedded in a simple group that can be embedded in a finitely presented group.
We prove the exact analogue for lattice-ordered groups:
\medskip

{\bf Theorem:}
{\it A finitely generated lattice-ordered group has soluble word problem if and only if it can be $\ell$-embedded in an $\ell$-simple lattice-ordered group that can be $\ell$-embedded in a finitely presented lattice-ordered group}.
\medskip

The proof uses permutation groups, a technique of Holland and McCleary, and the ideas used to prove the lattice-ordered group analogue of Higman's Embedding Theorem.
}
\end{abstract}
\bigskip

[{\bf Accepted and will appear in J. Group Theory in 2008.}]
\bigskip

-----------------------------------------


\newpage

\setcounter{theorem}{0}

\section{Introduction}

In 1974, W. W. Boone and G. Higman \cite{BH} proved:
\begin{theorem}
\label{bh}
{\rm (Boone-Higman \cite{BH})} A finitely generated group has soluble word problem if and only if it can be embedded in a simple group that can be embedded in a finitely presented group.
\end{theorem}

We prove the natural analogue for lattice-ordered groups:

\begin{theorem}
\label{lbh}
A finitely generated lattice-ordered group has soluble word problem if and only if it can be $\ell$-embedded in an $\ell$-simple lattice-ordered group that can be $\ell$-embedded in a finitely presented lattice-ordered group.
\end{theorem}

\noindent The proof of Theorem \ref{bh} was accomplished using $HNN$-extensions (spellings) and Higman's Embedding Theorem for groups:

\begin{theorem}
\label{h}
{\rm (Higman \cite{HI})} A finitely generated group can be embedded in a finitely presented group if and only if it can be defined by a recursively enumerable set of relations.
\end{theorem}

\noindent The difficult part of Theorem A was to find an algebraic condition equivalent to solubility of the word problem. The actual proof was relatively straightforward (given Theorem C).
In contrast (in the absence of spellings), our proof of Theorem \ref{lbh} uses a technique of Holland and McCleary \cite{HM} and the ideas of the proof of the lattice-ordered group analogue of Theorem \ref{h}:

\begin{theorem}
\label{lh}
{\rm (\cite{lH})} A finitely generated lattice-ordered group can be $\ell$-embedded in a finitely presented lattice-ordered group if and only if it can be defined by a recursively enumerable set of relations.
\end{theorem}

In Section 2, we give the basic background notation and results from previous papers in the subject and derive the easy half of the proof of Theorem \ref{lbh}.
In Section 3, we summarise the construction and formal proof of Theorem \ref{lh} from \cite{lH}, and in Section 4, we outline the permutation construction used there
and provide a modification.
In Section 5, we use this modification to consider the solubility of the word problem for a given recursively generated lattice-ordered group defined by a recursively enumerable set of relations. 
We use the results from Section 5 to deduce the harder half of Theorem \ref{lbh} in Section 6.
\medskip

To help the reader, I provide an outline of the proof.

We start with a finitely or recursively generated (lattice-ordered) group $G$ with soluble word problem.
Let $\{ (u_m,v_m):m\in \N\}$ be an enumeration of the pairs of (positive) non-identity elements of $G$.
Adjoin elements $\{ s_m:m\in \N\}$ and a modification of the relations $s_m^{-1}u_ms_m=v_m$ ($m\in \N$)\footnote{In the group case, the modification is only needed when the orders of $u_m$ and $v_m$ are different.  If the group is torsion-free, no modification is necessary.}.
Let $G^{\dagger}$ be the resulting countable (lattice-ordered) group.
To prove Theorem A, it suffices to show that $G^{\dagger}$ has soluble word problem (whence the argument is completed by continuing this construction inductively).
This is achieved by spelling/Britton extensions. In contrast, I have been unable to prove an analagous result directly in the lattice-ordered group case.
Instead, to prove Theorem B, I need to use the proof/construction of Theorem D and adjoin two extra elements $a_0,c_1$ and a finite set of relations.
To form $G^{\dagger}$, we also add the elements $\{ s_m:m\in \N\}$ and, \underline{\it inter alia}, relations $s_m^{-1}c_1^{-m}u_mc_1^ms_m=c_1^{-m}v_mc_1^{m}$ ($m\in \N$).
This may not be a free product with amalgamated $\ell$-subgroup as such a construction is not possible in general for lattice-ordered groups.
To complete this step of the proof of Theorem B, we need to do two things.
The first is to show that $G$ actually $\ell$-embeds in the recursively generated lattice-ordered group $G^{\dagger}$ which has a recursively enumerable set of defining relations. The second part of this main step is to prove that $G^{\dagger}$ has soluble word problem.  This is the crux of the proof. Since $G^{\dagger}$ is defined by a recursively enumerable set of relations, there is an algorithm to determine if a word $w$ in the alphabet of $G^{\dagger}$ is the identity. To find an {\it algorithm} that shows that a non-identity word is not the identity, we proceed by successively reducing the set of possible ``obstacle" words using wreath products, and then handling the remaining words by a technique due to Holland and McCleary \cite{HM}. The $\ell$-group $G^{\dagger}$ is specifically constructed for this technical part of the proof.

\section{Background and notation}

Throughout we will use $\N$ for the set of non-negative integers, $\Z_+$ for the set of positive integers, $\Q$ for the set of rational numbers and $\R$ for the set of real numbers. The only order on $\Q$ and $\R$ that we will consider will be the usual one.
\medskip

We assume that the reader has a minimal knowledge of recursive function theory (see \cite{Ro}).
\medskip

In any group $G$ we write $f*g$ for $g^{-1}fg$, and $[f,g]$ for $f^{-1}g^{-1}fg$. The former is often written $f^g$, though that would be less readable here where the expressions for $g$ are complicated.
Throughout, for any $m,n\in \Z_+$, we will write $f_1\dots f_m*g_1\dots g_n$ as a shorthand for $(f_1\dots f_m)*(g_1 \dots g_n)$.
\medskip

A {\it lattice-ordered group} is a group which is also a lattice that satisfies the identities
$x(y\wedge z)t= xyt\wedge xzt$ and $x(y\vee z)t= xyt\vee xzt$.
Throughout we write $x\leq y$ as a shorthand for $x\vee y=y$ or $x\wedge y=x$, 
and $\ell$-{\it group} as an abbreviation for lattice-ordered group.
A sublattice subgroup of an $\ell$-group is called an $\ell$-{\it subgroup}.

Lattice-ordered groups are torsion-free and $f\vee g=(f^{-1}\wedge g^{-1})^{-1}$.
Moreover, as lattices, they are distributive (\cite{G}, Lemma 2.3.5). Each element of $G$ can be written in the form $fg^{-1}$ where
$f,g\in G^+=\{ h\in G: h\geq 1\}$ --- see, {\it e.g.}, \cite{G}, Corollary 2.1.3, Lemma 2.3.2 \& Lemma 2.1.8. For each $g\in G$, let $|g|=g\vee g^{-1}$. 
Then $|g|\in G_+$ iff $g\neq1$, where $G_+=G^+\setminus \{ 1\}$.
Therefore, ($w_1=1~\&~\dots ~\&~w_n=1$) iff
$|w_1|\vee \dots \vee |w_n|=1$ [{\it ibid},  Lemma 2.3.8 \& Corollary 2.3.9].
Consequently, in the language of lattice-ordered groups (and in sharp contrast to group theory) any finite number of equalities can be replaced by a single equality.
Also, if $g\in G\setminus \{ 1\}$ and any element of $G_+$ is conjugate to $|g|$, then the normal subgroup of $G$ generated by $|g|$ is $G$; hence the normal $\ell$-subgroup of $G$ gnerated by $g$ is $G$ under this extra hypothesis.

We will write $f\perp g$ as a shorthand for $|f|\wedge |g|=1$ and say that $f$ and $g$ are {\it orthogonal}. As is well-known and easy to prove, $f \perp g$ implies $[f,g]=1$.

We will write $f\ll h$ if $f^m\leq h$ for all $m\in \Z$.

An $\ell$-{\it homomorphism} from one $\ell$-group to another is a group  and a lattice homomorphism. Kernels are precisely the normal $\ell$-subgroups that are convex (if $k_1,k_2$ belong to the kernel and $k_1\leq g\leq k_2$, then $g$ belongs to the kernel). They are called $\ell$-{\it ideals}. If the only $\ell$-ideals of an $\ell$-group are itself and $\{ 1\}$, then we say that the $\ell$-group is $\ell$-{\it simple}.
By the observation above, $G$ is $\ell$-simple if any two elements of $G_+$ are conjugate.

The free $\ell$-group on any set of generators exists by universal algebra. A finitely generated $\ell$-group is an $\ell$-homomorphic image of the free $\ell$-group on that finite number of generators.
If the kernel is finitely generated as an $\ell$-ideal, then we call the $\ell$-homomorphic image {\it finitely presented}; if the kernel is generated (as an $\ell$-ideal) by a recursively enumerable set of elements, then we say that the finitely generated $\ell$-homomorphic image has a {\it recursively enumerable set of defining relations}. We will write 
$$\langle Y~:~w_i(Y)=1~(i\in I)\rangle$$
for the quotient $F/K$ where $F$ is the free $\ell$-group on the generating set $Y$ and $K$ is the $\ell$-ideal generated (as an $\ell$-ideal) by
$\{ w_i(Y):i\in I\}$.
\medskip

The free $\ell$-group on a single generator is $\Z\oplus \Z$ ordered by:
$(m_1,m_2)\geq (0,0)$ iff $m_1,m_2\geq 0$; $(1,-1)$ is a generator since
$(1,-1)\vee (0,0)=(1,0)$.
\medskip 

We can already prove the easy half of Theorem \ref{lbh}.
The proof is identical to the group case (\cite{BH} or \cite{LS}, page 216).

\medskip

{\it Proof:} Suppose that a finitely generated $\ell$-group $G$ can be $\ell$-embedded in an $\ell$-simple $\ell$-group $S$ which can be $\ell$-embedded in a finitely presented $\ell$-group $H$.
Then $G$ has a recursively enumerable set of defining relations
since it can be $\ell$-embedded in a finitely presented $\ell$-group.
So, given $w\in G$, we can determine if $w=1$ in $G$.
To determine if $w\neq 1$ in $G$, let $g_1,\dots,g_n$ be the generators of $G$ and $g:=|g_1|\vee \dots \vee |g_n|$.
Let $\psi:G\rightarrow H$ be the $\ell$-embedding.
Let $H_w$ be the quotient of $H$ obtained by adjoining one extra relation, $w\psi=1$.
The $\ell$-ideal of $S$ generated by the image of $w$ must be all of $S$ if $w\neq 1$ in $G$, and is otherwise $\{ 1\}$. Thus $w\neq 1$ in $G$ 
if and only if $g\psi=1$ in $H_w$.
But $H_w$ is finitely presented, so we can determine if $g\psi =1$ in $H_w$.
The two algorithms together give the solubility of the word problem for $G$.
This proves the easy half of Theorem \ref{lbh}.
\pfend

The purpose of this article is to prove the converse (and so Theorem \ref{lbh}).
\medskip

In contrast to groups, the amalgamation property fails for $\ell$-groups: there are $\ell$-groups $G,H_1,H_2$ with $\ell$-embeddings $\sigma_j:G\rightarrow H_j$ ($j=1,2$) such that there is no $\ell$-group $L$ such that $H_j$ can be $\ell$-embedded in $L$ ($j=1,2$) so that the resulting diagram commutes (see \cite{P} or \cite{G}, Theorem 7.C). Hence $HNN$-extension techniques cannot be used (see \cite{G1}). Instead we use permutation group methods.
\medskip

Let
$(\Omega,\leq)$ be a totally ordered set. Then $Aut(\Omega,\leq)$ is an $\ell$-group when the group operation is composition and the lattice operations are just the pointwise supremum and infimum ($\alpha(f\vee g)=max\{ \alpha f, \alpha g\}$, {\it etc}.) 
There is an analogue of Cayley's Theorem for groups, namely the Cayley-Holland Theorem (\cite{G}, Theorem 7.A): 

\begin{theorem}
\label{ho}
{\rm (Holland \cite{HO})} Every lattice-ordered group can be $\ell$-embedded in $Aut(\Omega,\leq)$ for some totally ordered set $(\Omega,\leq)$; every countable lattice-ordered group can be be $\ell$-embedded in $Aut(\Q,\leq)$ and
hence in $Aut(\R,\leq)$.
\end{theorem}

We will write $A(\Omega)$ as a shorthand for $Aut(\Omega,\leq)$ when the total order on $\Omega$ is clear.
If $\De_1,\De_2\subseteq \Om$ are intervals, we will write $\De_1<\De_2$ if
$\de_1<\de_2$ for all $\de_j\in \De_j$ ($j=1,2$); we will write $\De_1\prec \De_2$
if $\De_1 < \Lambda < \De_2$ for some non-empty open interval $\La$ of $\Om$. 
\medskip

If $g\in A(\Omega)$, then the {\it support} of $g$, $supp(g)$, is the set
$\{\beta\in \Omega: \beta g\neq \beta\}$. 
\medskip

Since each real interval $(\alpha,\beta)$ is order-isomorphic to $(\R,\leq)$, we obtain:

\begin{corollary}
\label{supp}
Let $\alpha,\beta\in \R$ with $\alpha<\beta$. Then every countable $\ell$-group $G$ can be $\ell$-embedded in $A(\R)$ so that $supp(g)\subseteq (\alpha,\beta)$ for all $g\in G$.  
\end{corollary}

If $g\in A(\Omega)$ and $\alpha\in supp(g)$, then the convexification of the $g$-orbit of $\alpha$ is called the {\it interval of support of $g$ containing $\alpha$}; {\it i.e.}, the supporting interval of $g$ containing $\alpha$ is
$\{ \beta\in \Omega:~(\exists m,n\in \Z)( \alpha g^n\leq \beta\leq \alpha g^m)\}$.
So the support of an element is the disjoint union of its supporting intervals. 
The restriction of $g$ to one of its intervals of support is called a {\it bump} of $g$.
We will also call an element of $A(\Om)$ a {\it bump} if it has just one bump.
If $g$ is a bump, we write $\De_g$ for its unique supporting interval.
\medskip

By considering intervals of support, it is easy to establish the well-known fact:

\begin{proposition}
\label{conj}
For all $f,g\in A(\Omega)$, $supp(f*g)=supp(f)g$. Hence
if $f*g\perp f$ and $g\geq 1$, then $|f|\ll g$.
\end{proposition}

Let $\{ G_x : x\in X\}$ be a family of $\ell$-groups.  Then the full Cartesian product $C:=\prod \{G_x : x\in X\}$ is an $\ell$-group under the ordering
$$(g_x)_{x\in X}\in C^+ \;\;\;\hbox{iff}\;\;\; g_x\in G_x^+ \;\;\;\hbox{for all}\;\;\; x\in X.$$

\noindent We call $C$ the {\it cardinal product} of $\{G_x:x\in X\}$.
The restriction of this lattice-order to the direct sum $D:=\sum \{G_x:x\in X\}$ gives the {\it cardinal direct sum} (which is also an $\ell$-group). 
\medskip

Throughout, we will consider the restricted (small) wreath product (as opposed to the full Wreath product).
Let $(H,\Omega)$ be an $\ell$-permutation group; that is,
$H$ is an $\ell$-subgroup of $A(\Om)$.
We define the {\it wreath product} $W$ of an $\ell$-group $G$ and $(H,\Om)$,
written $G \wr (H,\Om)$, in the standard way: the base group, $B$, is
$\sum_{\alpha\in \Om} G_{\alpha}$, the direct sum (not full Cartesian product) of $\Om$ copies of $G$.
If $w:=(\{ g_{\alpha}\},h)\in W$, then $w\in W^+$ iff
$$h\in H^+\;\;\hbox{and}\;\; g_{\beta}\in G^+\;\hbox{for all}\; \beta\in \Om\;\;\hbox{with}\; \beta h=\beta.$$
As is standard, this makes $W$ into an $\ell$-group (see \cite{0}, Chapter 5) with the cardinal direct sum order on $B$.
 \medskip

We complete this section with two applications of the Cayley-Holland Theorem that we will need in the proof of Theorem B.

\subsection{Conjugacy}

Note that any conjugate of a strictly positive element of an $\ell$-group must be strictly positive. We first show that any two strictly positive bumps in $A(\R)$ of bounded support are conjugate and describe all conjugators.

\begin{lemma}
\label{1bump}
{\rm (\cite{HO})}
Let $f,g\in A(\R)_+$ be bumps with supports bounded above and below {\rm (}in $\R${\rm )}.
Let $\alpha,\beta\in \R$ be arbitrary with $\alpha\in supp(f)$ and $\beta\in supp(g)$.
Let $h_0:[\alpha,\alpha f]\rightarrow [\beta,\beta g]$ be any order-preserving bijection.
Then $h_0$ can be extended to an element $h\in A(\R)$ such that $h^{-1}fh=g$, and the restriction of $h\supseteq h_0$ to $supp(f)$ is uniquely determined.
\end{lemma}

{\it Outline of Proof:}
Let $m\in \Z$. Let 
$h_m:[\alpha f^m,\alpha f^{m+1}]\rightarrow [\beta g^m,\beta g^{m+1}]$ be given by
$h_m=f^{-m}h_0g^m$.
Let $h_*:supp(f)\rightarrow supp(g)$ be the union of these order-preserving bijections:
$h_*=\bigcup_{m\in \Z} h_m$.
Extend $h_*$ to an element $h \in A(\R)$ using arbitrary order-preserving bijections 
$(-\infty, inf\{ supp(f)\}]\rightarrow (-\infty, inf\{ supp(g)\}]$ and
$[sup\{ supp(f)\}, \infty) \rightarrow [sup\{ supp(g)\},\infty)$.
A simple calculation shows that $h^{-1}fh=g$ and that any $h\in A(\R)$
which conjugates $f$ to $g$ and extends $h_0$ must agree with $h_m$ on
$[\alpha f^m,\alpha f^{m+1}]$ ($m\in \Z$) and so extends $h_*$.
\pfend

The key here is that $\alpha$ and $\beta$ are arbitrary in the supports of $f$ and $g$ respectively, and so is the order-preserving bijection $h_0$ from $[\alpha,\alpha f]$ to $[\beta,\beta g]$.
\medskip

For any $f\in A(\Om)_+$, let ${\mathcal B}_f$ be the set of bumps of $f$.
Let ${\mathcal S}_f:=\{ \De_{f'}:f'\in {\mathcal B}_f\}$ be the set of supports of bumps of $f$. Then ${\mathcal S}_f$ inherits the interval order from $\Om$; {\it i.e.}, $\De_{f_1}<\De_{f_2}$ if this holds in $\Om$ 
($f_1,f_2\in {\mathcal B}_f$).
Let $({\mathcal I}_f,<)$ be the set of intervals of $\Om$ maximal with respect to being disjoint from $supp(f)$, equipped with the inherited interval order from $\Om$.
Let $\Lambda_f:={\mathcal S}_f \cup {\mathcal I}_f$ with the inherited interval order from $\Om$.
\medskip

The proof of Lemma \ref{1bump} extends to show

\begin{lemma}
\label{conjbumps}
{\rm (\cite{HO})}
Let $f,g\in A(\R)_+$. 
Suppose that there is is an order-preserving bijection 
$\varphi:(\Lambda_f,<) \rightarrow (\Lambda_g,<)$ that restricts to a bijection between ${\mathcal S}_f$ and ${\mathcal S}_g$.
For each $f_j\in {\mathcal B}_f$, let $g_j \in {\mathcal B}_g$ be such that
$\De_{f_j}\varphi =\De_{g_j}$ ($j\in J$).
Let $\alpha_j\in supp(f_j)$ and $\beta_j \in supp(g_j)$ be arbitrary and 
$h_{j,0}:[\alpha_j,\alpha_j f]\rightarrow [\beta_j, \beta_j g]$ be an arbitrary order-preserving bijection.
Then there is $h\in A(\R)$, uniquely defined on $supp(f)$, that extends all $h_{j,0}$
and conjugates all $f_j$ to $g_j$ ($j\in J$), and so conjugates $f$ to $g$.
\end{lemma}

{\bf Remark}. We will later use such freedom of choice for $\varphi$ to show that certain $\ell$-group expressions cannot be the identity.
\medskip

Observe that if $f,g\in A(\R)_+$ and $f$ has a single bounded bump but $g$ has more than one bump, then $f$ and $g$ are not conjugate in $A(\R)$.
Nonetheless, although the amalgamation property fails for $\ell$-groups, Keith Pierce \cite{P} was able to use the Cayley-Holland Theorem to prove:

\begin{theorem}
\label{p}
{\rm (K. R. Pierce \cite{P})}
Every $\ell$-group can be $\ell$-embedded in one in which any two strictly positive lements are conjugate.
\end{theorem}

In order to prove Theorem \ref{lbh}, we will need some of the ideas of the proof of Theorem \ref{p}. We therefore provide a very brief sketch of Pierce's proof here.
\medskip

First observe that it suffices to prove that every $\ell$-group $G$ can be $\ell$-embedded in an $\ell$-group $H$ in which the images of any two strictly positive elements of $G$ are conjugate. For if $G^{\dagger}$ is the $\ell$-subgroup of
$H$ generated by the image of $G$ and the conjugators in $H$, define $G(0)=G$ and $G(m+1)=G(m)^{\dagger}$ ($m\in \N$).
Let $\hat G:=\bigcup_{m\in \N} G(m)$. Then $G$ can be $\ell$-embedded in $\hat G$ and any two elements of  $\hat G_+$ are conjugate in $\hat G$.
\medskip

By $\ell$-embedding $G$ diagonally into the cardinal product $\prod \{ G:n\in \Z\}$ if necessary, we may assume that for any $f,g\in G_+$, no supporting interval $\De$ of $g$ is greater than all supporting intervals of $f$, nor less than all supporting intervals of $f$.

Next, by a modification of the Cayley-Holland Theorem due to Weinberg \cite{W1}, we may assume that the $\ell$-group $G$ is contained in $B(T)$, the $\ell$-group of all order-preserving bijections of bounded support for some totally order set $(T,\leq)$ in which, for any $\sigma_j<\tau_j$ in $T$ ($j=1,2$), there is $h\in B(T)$ such that $\sigma_1 h=\sigma_2$ and $\tau_1 h=\tau_2$. We may therefore assume that the $\ell$-group $G$ is $B(T)$ for such a totally ordered set $(T,\leq)$ (see {\it op. cit.} or \cite{0}, Corollary 2.L.).
\medskip

As noted above, we need to be able to $\ell$-embed $G$ in some
$A(\Xi)$ so that there is bijection between $\Lambda_f$ and $\Lambda_g$ in $\Xi$ for any $f,g$ images of elements of $B(T)_+$. We have one further complication which did not arise in the case of $A(\R)$; the endpoints of a bump in the Dedekind completion $\bar \Xi$ of $\Xi$ may or may not belong to $\Xi$, or to the same orbit of $A(\Xi)$ in the natural action of $A(\Xi)$ on $\bar \Xi$. So we will also need to construct the $\ell$-embedding so that the lower endpoints of corresponding bumps of $f$ and $g$ must belong to the same $A(\Xi)$ orbit, and ditto for upper endpoints of corresponding bumps.
This is achieved via transfinite induction assuming that 
$2^{|B(T)|}=|B(T)|^+$, the successor cardinal of the cardinality of $B(T)$.

At even successor stages, one employs the orbit Wreath product;
at odd successor stages, one adjoins certain cuts of the Dedekind completion of the previous totally ordered set and extends the order-preserving permutations uniquely; and at limit stages, one takes unions.
For the details, see \cite{P} or \cite{0}, pp. 194-205.
At stage $|B(T)|$, we obtain an $\ell$-permutation group $(H,\Om)$ with $T\subseteq \Om$ and $H\subseteq B(\Om)$.
We identify $G$ with its image in $H\subseteq B(\Om)$. 
Now the set ${\mathcal I}_g$ of fixed point intervals of any $g\in G_+$ has a greatest and least element.
Let ${\mathcal I}^-_g$ be the totally ordered subset of ${\mathcal I}_g$ obtained by removing these two intervals,
and $\Lambda^-_g={\mathcal I}^-_g \;\cup \; {\mathcal S}_g$.
Then, in $\Lambda^-_g$, the pair $({\mathcal I}^-_g, \; {\mathcal S}_g)$ forms a $|B(T)|$-set of type $2$.
For any $f,g\in G_+$, there is $\phi:\Lambda_f\rightarrow \Lambda_g$ such that for all $x,y\in G_+$ with $x\not\perp f$ and $y\perp f$
\medskip

(i) 
for uncountably many intervals $\Delta_{f'}\in {\mathcal S}_f$ with $\De_{f'}\cap \De_{x'}\neq \emptyset$ for some $\De_{x'}\in {\mathcal S}_x$, there is  $\De_{y'}\in {\mathcal S}_y$
such that $$\De_{f'}\phi \prec \De_{y'} \prec \De_{f'},\; \hbox{and}$$

(ii) for uncountably many intervals $\Delta_{f'}\in {\mathcal S}_f$ with $\De_{f'}\cap \De_{x'}\neq \emptyset$ for some $\De_{x'}\in {\mathcal S}_x$, there is  $\De_{y'}\in {\mathcal S}_y$
such that $$\De_{f'}\phi \succ \De_{y'} \succ \De_{f'}.$$

 (iii) If the supports of $f$ and $g$ are not disjoint, then we can also require that  
$\Delta_{f'} \phi\;\cap \;\Delta_{f'}\neq \emptyset$
for uncountably many intervals $\Delta_{f'}\in {\mathcal S}_f$.
\medskip

We can use a natural extension of Lemma \ref{conjbumps} to obtain
$h\in H$ conjugating $f$ to $g$ such that  for all $x,y\in G_+$ with $x\not\perp f$ and $y\perp f$
\medskip

(I) for uncountably many intervals $\Delta_{f'}\in {\mathcal S}_f$ with $\De_{f'}\cap \De_{x'}\neq \emptyset$ for some $\De_{x'}\in {\mathcal S}_x$, there is  $\De_{y'}\in {\mathcal S}_y$
such that $$\De_{f'}h \prec \De_{y'} \prec \De_{f'},\; \hbox{and}$$

(II) for uncountably many intervals $\Delta_{f'}\in {\mathcal S}_f$ with $\De_{f'}\cap \De_{x'}\neq \emptyset$ for some $\De_{x'}\in {\mathcal S}_x$, there is  $\De_{y'}\in {\mathcal S}_y$
such that $$\De_{f'}h \succ \De_{y'} \succ \De_{f'}.$$

(III) if the supports of $f$ and $g$ are not disjoint, then 
$\Delta_{f'} h \;\cap \;\Delta_{f'}\neq \emptyset$
for uncountably many intervals $\Delta_{f'}\in {\mathcal S}_f$.
\medskip

For more details, see \cite{P} or \cite{0}, pp.194-205.
\medskip

Although it is not explicit, the proof yields further information.

\noindent Suppose that for some $c\in G_+$ we have, for each supporting interval $\De$ of $c$ in $T$, there is
$\alpha_{\De}\in T$ such that 
$$\alpha_{\De} \prec (supp(f) \;\cup\; supp(g)) \; \cap \; \De \prec \alpha_{\De} c.$$
In the extension of $(G,T)$ to $(H,\Om)$ we can ensure that, for every supporting interval $\De$ of $c$ in $\Om$, 
$\;\;\;\alpha_\De \prec (supp(f)\;\cup\; supp(g))\;\cap \; \De \prec \alpha_\De c$.

\noindent Indeed, Pierce's construction ensures that $f$ and $g$ have uncountably many supporting intervals contained in $(\alpha_{\De}, \alpha_{\De})$.
We then have, for each $m\in \N$, a map 
$$\phi_m:(\Lambda_{f*c^m},<) \rightarrow (\Lambda_{g*c^m},<)$$ 
\noindent with
${\mathcal S}_{f*c^m}\phi_m={\mathcal S}_{g*c^m}$ and
${\mathcal I}^-_{f*c^m}\phi_m={\mathcal I}^-_{g*c^m}$, so that $\phi_m$ has the properties (i) -- (iii) above with $f*c^m$ in place of $f$ and $g*c^m$ in place of $g$.
We can 
further ensure that $\phi_m$ induces the identity off 
$\; \bigcup \{ (\alpha_{\De}c^m,\alpha_{\De} c^{m+1}):m\in \N,\; \De\in {\mathcal S}_c\}.$

\noindent
As above, there is a resulting pairwise orthogonal set $\{h_m\in H: m\in \N\}$ 
with $\Delta h_m=\De \phi_m$ ($\De\in {\mathcal S}_{f*c^m}$ or
$\De\in {\mathcal I}_{f*c^m}$) 
such that properties (I) -- (III) above hold with $f*c^m$ in place of $f$ and $g*c^m$ in place of $g$; so 
$$f*c^mh_m =g*c^m \;\;\hbox{(}m\in \N\hbox{)}\;\;\hbox{and}\;\;
f*c^{m'}h_m=f*c^{m'}\;\;\hbox{(}m,m'\in \N,\;m'\neq m\hbox{)}.$$

E. C. Weinberg (\cite{W2} and \cite{W3}) has shown how to remove all dependence on any form of the Generalised Continuum Hypothesis.
This is achieved by using Harzheim's minimal $\eta_{\kappa}$-sets instead of $\kappa$-sets where $\kappa=|B(T)|$ (see \cite{HA}). 
The proof proceeds exactly as before with this minor modification at stages and provides conjugators with the same properties.
\bigskip

\subsection{The word problem for free $\ell$-groups}

Another application of the Cayley-Holland Theorem was provided independently by Kopytov and McCleary. They proved that the free lattice-ordered group on a finite number of generators has a faithful highly transitive representation (\cite{K}, \cite{M} or \cite{G}, Theorem 8.D). Indeed,

\begin{proposition}
\label{HM}
{\rm \cite{HM}}. Given any order-preserving isomorphisms $z_j$ with domain and range {\it finite} subsets of $\R$ {\rm(}$j=1,\dots,n${\rm)}, these maps can be extended to elements $y_j\in A(\R)$ {\rm(}$j=1,\dots,n${\rm)} so that the $\ell$-subgroup of $A(\R)$ generated by $\{ y_1,\dots,y_n\}$ is the free $\ell$-group $F$ on $\{ y_1,\dots,y_n\}$.
\end{proposition}

Holland and McCleary applied  this to prove ({\it op. cit.})

\begin{theorem}
\label{free}
{\rm \cite{HM}}
For any positive integer $n$, the free lattice-ordered group on $n$ free generators has soluble word problem.
\end{theorem}

The idea of the proof is as follows.

First consider a single group term $w(y_1,\dots,y_n)$, say $w:=y_{j_1}^{\epsilon_1}\dots y_{j_k}^{\epsilon_k}$, where
$j_1,\dots , j_k\in \{ 1,\dots,n\}$ and $\epsilon_1,\dots,\epsilon_k\in \{\pm 1\}$.
We draw two diagrams, one with $0 y_{j_1}>0$, the other with
$0 y_{j_1}<0$.

From each of these diagrams we construct three new diagrams
if $j_2\neq j_1$. 
For the first diagram ($0y_{j_1}>0$), we make the following modification. If $\epsilon_1=1$, we draw three diagrams, the first with 
$$0y_{j_1}^{\epsilon_1}y_{j_2}^{\epsilon_2}>0y_{j_1}^{\epsilon_1}>0,$$
the second with $$0y_{j_1}^{\epsilon_1}>0y_{j_1}^{\epsilon_1}y_{j_2}^{\epsilon_2}>0,$$
and the third with 
$$0y_{j_1}^{\epsilon_1}>0>0y_{j_1}^{\epsilon_1}y_{j_2}^{\epsilon_2};$$
on the other hand, if $\epsilon_1=-1$, we construct three diagrams:
in the first, we have
$$0y_{j_1}^{\epsilon_1}y_{j_2}^{\epsilon_2}>0>0y_{j,1}^{\epsilon_1},$$
in the second 
$$0>0y_{j_1}^{\epsilon_1}y_{j_2}^{\epsilon_2}>0y_{j,1}^{\epsilon_1},$$
and in the third
$$0>0y_{j_1}^{\epsilon_1}>0y_{j_1}^{\epsilon_1}y_{j_2}^{\epsilon_2}.$$

If $0y_{j_1}>0$ and $j_1=j_2$, then we construct a single diagram with
$$0y_{j_1}^{\epsilon_1}y_{j_2}^{\epsilon_2}>0y_{j_1}^{\epsilon_1}>0
\;\;\; \hbox{if}\;\; \epsilon_1=\epsilon_2=1;$$ a single diagram with
$$0y_{j_1}^{\epsilon_1}y_{j_2}^{\epsilon_2}<0y_{j_1}^{\epsilon_1}<0\;\;\;
\hbox{if}\;\;\epsilon_1=\epsilon_2=-1;$$ 
a single diagram with $$0y_{j_1}^{\epsilon_1}y_{j_2}^{\epsilon_2}=0<0y_{j_1}^{\epsilon_1}
\;\;\;\hbox{if}\; \epsilon_1=1\;\hbox{and}\;\epsilon_2=-1;$$
and a single diagram with $$0y_{j_1}^{\epsilon_1}y_{j_2}^{\epsilon_2}=0>0y_{j_1}^{\epsilon_1}
\;\;\;\hbox{if}\; \epsilon_1=-1\;\hbox{and}\;\epsilon_2=1.$$

Similarly, we construct diagrams from the second case ($0y_{j_1}<0$).
We proceed with the spelling ensuring only that when we consider $y_{j_i}^{\epsilon_i}$, the element $y_{j_i}$ and its inverse respect all the inequalities declared previously involving $y_{\ell}$ where $j_i=\ell$.

By Proposition \ref{HM}, if in {\it all} possible resulting legitimate diagrams we have $0w=0$, then $w=1$ in $F$;
if in {\it some} resulting legitimate diagram we get $0w\neq 0$, then $w\neq 1$ in $F$ by the same proposition.

This completes the solubility of the group word problem in $F$.
\medskip

For a general $\ell$-group word $w(y_1,\dots,y_n)$, enumerate the group words used to constitute $$w:=\bigvee_{i=1}^k \bigwedge_{j=1}^{r_i} w_{i,j};$$ {\it i.e.}, $w_{1,1},\dots,w_{1,r_1},w_{2,1},\dots,w_{k,r_k}.$
Form all possible legitimate diagrams as above for $w_{1,1}$.
For each of these diagrams, do the same for $w_{1,2}$ subject only that all inequalities that occurred in  that diagram for $w_{1,1}$ are respected in the diagrams for $0w_{1,2}$.
For each of the resulting diagrams, do the same for $w_{1,3}$, {\it etc.}
Then $w=1$ in $F$ if $0w=0$ in {\it all} resulting diagrams;
and $w\neq 1$ in $F$ if $0w\neq 0$ in {\it some} resulting diagram.
\pfend

We will use the idea of this proof in the last part of the proof of Theorem \ref{lbh}.

\section{Summary of the proof of Theorem \ref{lh}.}

The proof extends the ideas in \cite{GG}.
\medskip

Let $H$ be an $\ell$-group that has generators $\{ y_n:n\in \Z_+\}$ and is defined by a recursively enumerable set of relations. 
Then there is an algorithm that constructs a $2$-generator $\ell$-group $\bar H$ and an explicit $\ell$-embedding of $H$ into $\bar H$ such that (the image of) every element of $H$  is equal to a group term in the generators of $\bar H$ and $\bar H$ is definable by a recursively enumerable set of $\ell$-group words; the defining relations for $\bar H$ are group terms or finite meets of group terms and are explicitly obtainable from the defining relations of $H$ (see the proof of Theorem E in \cite{lH}, Section 6). Moreover, the proof in \cite{lH} shows that this set of defining relations for $\bar H$ is recursive if the set of defining relations for $H$ is, and $\bar H$ has soluble word problem whenever $H$ does.
\medskip

{\it We may therefore assume that $H$ is finitely generated with a recursively enumerable set of defining relations, each of which is equal to a group term or a finite meet of group terms; moreover, every element of $H$ is a group term in the generators. }
\medskip

In \cite{lH}, Section 3, we called these $\ell$-group words {\it meet strings} and gave an explicit recursive G\"{o}del numbering for the set of all meet strings occurring in the free $\ell$-group on the $n$ free generators $y_1,\dots,y_n$: for each meet string $w(y_1,\dots,y_n)$, we defined the G\"{o}del number $\ga(w)$ of $w$. 
Not all natural numbers were G\"{o}del numbers of meet strings. We rectified matters by providing an explicit recursive pseudo-G\"{o}del numbering for the set of all meet strings occurring in the free $\ell$-group on $y_1,\dots,y_n$ (\cite{lH}, Section 3);
each natural number was a pseudo-G\"{o}del number of a unique meet string and each non-empty meet string had an infinite recursive set of pseudo-G\"{o}del numbers. 
\medskip

For the $\ell$-group $H$ generated by $y_1,\dots,y_n$ and defined by a recursively enumerable set of meet string relations, let $X$ be the set of {\it all} pseudo-G\"{o}del numbers of {\it all} the meet strings in $y_1,\dots, y_n$ that hold in $H$. In \cite{lH}, Section 5.1, we constructed from $H$ (and $X$) a finitely presented $\ell$-group $L(X)$, and provided an explicit map $\varphi$ of $H$ into $L(X)$. In Section 5.2 of \cite{lH}, we proved that $\varphi$ was a well-defined $\ell$-homomorphism. 
\medskip

Crucially for our needs, there were generators $a_0,c_1\in L(X)_+$ such that, in $L(X)$, we had 

\begin{equation}
\label{00}
a_0*c_1^m\perp a_0,\;\;\;\;\;\hbox{for all}\;\; m\in \Z_+, 
\end{equation}
\noindent and for all distinct $m,m'\in \Z$, 
\begin{equation}
\label{pwd}
(y_j\varphi)*a_0^m \perp (y_k\varphi)*a_0^{m'}\;\;\;\;\;\&\;\;\;\;\;
(y_j\varphi)*c_1^m \perp (y_k\varphi)*c_1^{m'},
\end{equation}

\noindent where ($j,k\in \{ 1,\dots, n\}$).
\medskip

To show that $\varphi$ is injective, we used the Cayley-Holland Theorem to get a representation $\widehat{L(X)}$ of $L(X)$ which was faithful on $H\varphi$.
We briefly describe this in the next section.

\section{The permutation representation in \cite{lH}.}

In \cite{lH}, Section 5.3, we constructed order-preserving permutations of $\R$ that satisfied all the defining relations of $L(X)$. That is, we constructed an $\ell$-subgroup of $A(\R)$ that was an $\ell$-homomorphic image $\widehat{L(X)}$ of $L(X)$.
We proved that this permutation representation of $L(X)$ lead to a faithful representation for $H$.
By arrow chasing, it followed that the well-defined $\ell$-homomorphism $\varphi:H\rightarrow L(X)$ was injective. This proved Theorem \ref{lh}.
As in \cite{lH}, we identify $H$ with its image in $L(X)$; {\it i.e.,}
we take $\varphi$ to be the identity.
\medskip

In the presentation of $L(X)$, we had a generator $y$ such that for each $\ell$-group term $w(y_1,\dots,y_n)$, $y$ and $wy$ were (explicitly) conjugate in $L(X)$.
No attempt was made to try to conjugate $\ell$-group terms in  $y_1,\dots,y_n$ to each other in $L(X)$ (if they were strictly positive in $H$).
It was unnecessary in \cite{lH}. However, we will need to do so in this article to prove Theorem \ref{lbh}. 
We will add extra generators and relations to those of $L(X)$ to ensure that any two strictly positive elements of $H$ are conjugate in the new $\ell$-group (which will be countable). This is easy to achieve by Theorem \ref{p}.
We will want the induced $\ell$-homomorphism of $H$ into the constructed $\ell$-group to be injective. This will require modifying the permutation representation in \cite{lH}.
\medskip

Let $\Om$ be a minimal $\eta_1$-set. Instead of representing $L$ in $A(\R)$, we represent it in $A(\Om)$ so that the set of bumps and fixed point intervals of every $w\in H_+$ form a minimal $\eta_1$-set of type $2$.
In particular, as noted in Section 2.1, for each $w\in H$ and supporting interval $\De$ of $\hat c_1$, there is $\alpha_{\De}\in \Om$ such that
$$\alpha_{\De} \prec \De \;\cap\; supp(\hat w) \prec \alpha_{\De} \hat a_0 \prec \alpha_{\De} \hat c_1,$$
\noindent and there are uncountably many bumps of $\hat w$ in $(\alpha_{\De},\alpha_{\De}\hat a_0)$.
 This provides an $\ell$-embedding of $H$ into $B(\Om)$.

\section{Soluble word problem.}

Our aim in this section is to prove

\begin{proposition}
\label{prop}
Let $G$ be a recursively generated $\ell$-group defined by a recursively enumerable set of relations. Suppose that $G$ has soluble word problem.
Then $G$ can be $\ell$-embedded in a recursively generated $\ell$-group $G^{\dagger}$ with soluble word problem in which any two strictly positive elements of $G$ are conjugate. 
\end{proposition}

Throughout this section, let $G=F/K$ be a fixed recursively generated $\ell$-group with soluble word problem.
That is, $F$ is a free $\ell$-group on a recursive set of free generators (say, $\{ y_n:n\in \N\}$) and $K$ is an $\ell$-ideal such that the set of $\ell$-group terms in $F$ which belong to $K$ is recursive. 
\medskip

Now for each $w\in F$, we have $Kw\in G_+$ if and only if 
($w\not\in K$ but $w\wedge 1\in K$).
Since $K$ is recursive, we have an algorithm to determine whether or not $Kw\in G_+$. Hence

\begin{lemma}
\label{+}
If a recursively generated $\ell$-group $G$ has soluble word problem, then the strict positivity problem for $G$ is also soluble.
\end{lemma}

The following fact is folk-lore; a proof is included only because I have been unable to find one in the literature.

\begin{lemma}
\label{sum}
Let $X$ be a recursive set and $\{G_x:x\in X\}$ be a family of recursively generated $\ell$-groups, each with soluble word problem.  Then the cardinal sum $D$ of $\{G_x:x\in X\}$ has soluble word problem.
\end{lemma}

{\it Proof:} We assume that $G_x$ and $G_{x^\prime}$ share no common symbol except $1$.
Let $\{g_{x,m}:m\in\N\}$ generate $G_x$ and ${\mathcal R}_x$ be the recursive set of relations for $G_x\;\;(x\in X)$.  
Then $\{g_{x,m}:x\in X, m\in\N\}$ is a recursive set of generators for $D$.  The defining relations for $D$ are $\bigcup\{{\mathcal R}_x:x\in X\}$ together with $g_{x,m}\perp g_{x^\prime,m^\prime}\;\;(m,m^\prime\in\N;\; x,x^\prime\in X,\; x\neq x^\prime)$.
Thus $[g_{x,m},g_{x^\prime,m^\prime}]=1$ for all $m,m^\prime\in\N$ and distinct $x,x^\prime\in X$.

Each group word in $D$ has the form $w_{x_1}\dots w_{x_k}$ for some $x_1,\dots,x_k\in X$ distinct, where $w_{x_i}\in G_{x_i}$ are group terms $(i=1,\dots,k)$.  Thus any $\ell$-group word $w$ in the alphabet of $D$ has form $u_{x_1}\dots u_{x_k}$ for some $x_1,\dots,x_k\in X$ distinct, where $u_{x_i}\in G_{x_i}$ are $\ell$-group terms $(i=1,\dots,k)$.  For example, $g_{x,1}g_{x^\prime,2}^{-1}\wedge g_{x,3} = (g_{x,1}\wedge g_{x,3})(g_{x^\prime,2}^{-1}\wedge 1)$.  Since $D$ is recursively generated and defined by a recursively enumerable set of relations, there is an algorithm to determine if $w=1$ in $D$.  To determine if $w\neq 1$ in $D$, we need only check the equivalent fact that $u_{x_i}\neq 1$ in $G_{x_i}$ for some $i\in\{1,\dots,k\}$.  Since $w$ provides $x_1,\dots,x_k$ and each $G_{x_i}$ has soluble word problem, we can determine algorithmically if $u_{x_1}\neq 1 \;\hbox{or}\;
u_{x_2}\neq 1\;\hbox{or}\;\dots\;\hbox{or}\;u_{x_k}\neq 1$.  If all of these fail, then $w=1$ in $D$; if at least one of them holds, then $w\neq 1$ in $D$.  Thus $D$ has soluble word problem.\pfend

We need another well-known fact:

\begin{lemma}
\label{wr}
Let $G$ be a recursively generated $\ell$-group with soluble word problem and $c>1$ be a new symbol. Let $W$ be the $\ell$-group wreath product $G \wr (\langle c\rangle,\Z)$.
Then $W$ is a recursively generated $\ell$-group with soluble word problem.
\end{lemma}

{\it Proof:} Let $\{g_n:n\in\N\}$ generate $G$.  Then $\{c\}\cup\{g_n:n\in\N\}$ generates $W$ and $W$ is defined by the defining relations of $G$ together with
$$c\wedge 1=1,\;\;\;\;\;g_n * c^m \perp g_{n^\prime}\;\;\;(n,n^\prime\in\N; \;m\in\Z\setminus\{0\}).$$

\noindent Hence $W$ is recursively generated and defined by a recursively enumerable set of relations.  Let $w=\bigvee_I\bigwedge_J w_{i,j}$ be an $\ell$-group term in the generators of $W$ with each $w_{i,j}$ being a group term therein.  We have an algorithm to determine if $w=1$ (since $W$ is recursively generated and defined by a recursively enumerable set of relations).  To determine if $w\neq 1$ in $W$, put every element of $G$ occurring in $w$ equal to 1.  Let $w^\prime$ be the result.  Then $w^\prime\in\langle c\rangle\cong\Z$ so we can determine whether or not $w^\prime=1$ in $\langle c \rangle$.  If $w^\prime \neq 1$ in $\langle c\rangle$, then as $\langle c\rangle$ is the $\ell$-homomorphic image of $W$ with kernel the base group $B$, we have that $w\neq 1$ in $W$.  So assume that $w^\prime=1$ in $\langle c\rangle$, {\it i.e.}, $w\in B$.

If for some $i_0\in I$, there is $j_0\in J$ such that the sum of the exponents of $c$ appearing in $w_{i_0,j_0}$ (called the {\it weight} of $c$ in $w_{i_0,j_0}$) is negative, then $\bigwedge_J w_{i_0,j} <1$ in $W$.  Since $w'=1$ in $W/B$, we have that 
$\bigvee_{I\setminus \{ i_0\}}\bigwedge_J w_{i,j}=
\bigvee_I\bigwedge_J w_{i,j}=1$ in $W/B$.
Hence we may assume that the weight of $c$ in each $w_{i,j}$ is non-negative.
For each $i\in I$, let 
$$J_{i,0}:=\{ j\in J:\hbox{the weight of}\;c\;\hbox{in}\;w_{i,j} \;\hbox{is}\; 0\}.$$
Since $w'=1$ in $W/B$, we have $J_{i,0}\neq \emptyset$ for all $i\in I$.
If for some $i_1\in I$ we have $J_{i_1,0}\neq J$, then $\bigwedge_J w_{i_1, j}=\bigwedge_{J\setminus J_{i_1,0}} w_{i_1,j}$.
Hence we may assume that $J_{i,0}=J$ for all $i\in I$.
That is, the weight of $c$ in each $w_{i,j}$ is $0$.
So each $w_{i,j}$ is an element of $B$.
But $B=\sum \{ c^{-m}Gc^m:m\in \Z\}$. By Lemma \ref{sum}, $B$ has soluble word problem. We can therefore determine whether or not $w=1$ in $B$ and thus solve the word problem for $W$.\pfend

If $X$ and $Y$ are totally ordered sets, define $X \lex Y$ to be the set $X\times Y$ totally ordered by: $(x,y)<(x^\prime,y^\prime)$ if either $(y<y^\prime\;\hbox{in}\;Y)$ or $(y=y^\prime \;\hbox{in}\; Y\;\&\;x<x^\prime\;\hbox{in}\;X)$.
\medskip

We generalise Lemma \ref{wr} slightly.

\begin{lemma}
\label{Gwr}
Let $X=\Z\lex\Z$ and $A=\langle a\rangle\wr(\langle c\rangle,\Z)$ viewed as an $\ell$-subgroup of $A(X)$ in the natural way.
If $G$ is a recursively generated $\ell$-group with soluble word problem, then
$G \wr (A,X)$ has soluble word problem.
\end{lemma}

{\it Proof:}
This follows from two applications of Lemma \ref{wr}, since $G\wr (A,X)\cong (G\wr (\langle a\rangle,\Z))\wr (\langle c\rangle,\Z)$.\pfend

Since $G$ is recursively generated and has soluble word problem, it
can be defined by a recursively enumerable set of relations and
so can be $\ell$-embedded in $L$ as in \cite{lH}.
This construction and $\ell$-embedding were effective (as noted in Sections 3 and 4).
We identified $G$ with its $\ell$-isomorphic image in $L$.
The $\ell$-subgroup of $L$ generated by $G_{\flat}:=G\cup \{ a_0,c_1\}$ is 
$\ell$-isomorphic to $G\wr (\langle a_0\rangle,\Z) \wr (\langle c_1\rangle,\Z)$.
Hence, by Lemma \ref{Gwr}, 

\begin{lemma}
\label{flat}
If the recursively generated $\ell$-group $G$ has soluble word problem, then $G_{\flat}$ has soluble word problem.
\end{lemma}

We will also introduce an infinite set of conjugators, $\{ s_m:m\in \N\}$.
By Lemma \ref{sum},

\begin{lemma}
\label{S}
Let $\{ s_m:m\in \N\}$ be a set of new symbols and for each $m\in \N$, let $S_m$ be the free $\ell$-group on the single free generator $s_m$.
Let $S=\sum_{m\in \N} S_m$ be the (abelian) $\ell$-group with the cardinal ordering.
Then $S$ is recursively generated, recursively defined and has soluble word problem.
\end{lemma}

{\bf Caution}: The construction and consequent proof of Proposition \ref{prop} is complicated by the failure of the amalgamation property. We must ensure that not too much ``collapses" so that $G$ is still $\ell$-embeddable in the resulting $\ell$-group.
\medskip

As noted in Section 3, we may assume that every element of $G$ can be written as a group word in the generators. 
By Lemma \ref{+}, there is a recursive enumeration $(u_0,v_0),(u_1,v_1),(u_2,v_2),\dots$ of all pairs of elements of $G_+$.
Let $G^{\dagger}$ be the $\ell$-homomorphic image of the $\ell$-group free product of $G_{\flat}$ and $S$ obtained by adjoining the recursive set of extra relations 
\begin{equation}
\label{e1}
|s_m|^k \leq a_0*c_1^m\;\;\;\;\;(m,k\in \N),
\end{equation}

\begin{equation}
\label{e3}
s_m*c_1^{-m}\perp g_n*a_0^k\;\;\;\;\;\hbox{(}m,n\in \N,\; k\in \Z\setminus \{ 0\}\hbox{)}, 
\end{equation}

\begin{equation}
\label{e5}
s_m*c_1^{-m}\perp s_{m'}*c_1^{-m'}a_0^k\;\;\;\;\;\hbox{(}m,m'\in \N,\; k\in \Z\setminus \{ 0\}\hbox{)}, 
\end{equation}

\begin{equation}
\label{e4}
u_m*c_1^ms_m=v_m*c_1^m\;\;\;\;\;(m\in \N),
\end{equation}
\medskip

\noindent where $\{ g_n:n\in \N\}$ generates $G$.  Note that, since $\{ a_0 * c_1^m: m\in \N\}$ is a pairwise orthogonal set of elements of $G_{\flat}$ (and hence of $G^{\dagger}$), we have

\begin{equation}
\label{e2}
 s_m\perp a_0*c_1^{m'}\;\;\;\; \hbox{(}m,m'\in \N,\;\;\;m'\neq m\hbox{)}, 
\end{equation}

\noindent Indeed, the relations imply that for all $u\in G\cup S$, we have
\begin{equation}
\label{orth}
(|u| \wedge (a_0 * c_1^m))^k \leq a_0 * c_1^m \;\;\;\; \hbox{for all} \; m,k\in\Z.
\end{equation}

For the remainder of this section, let $S^\prime$ be the $\ell$-subgroup of $G^{\dagger}$ generated by $\{ s_m * c_1^{-m} : m \in \N\}$.
\medskip

As noted at the end of Section 4, (using the notation of Section 2.1) we have set up the relations (in $B(\Om)$) to get $\ell$-homomorphic images of $G_{\flat}$ and $S$ in $B(\Om)$
such that, for all $m\in \N$,
 $$\hat u_m*\hat c_1^m \hat s_m=\hat v_m*\hat c_1^m,\;\;\hbox{and}$$
$$\alpha_{\De}\hat c_1^m \prec \De \;\cap\; supp(\hat s_m)  \prec \alpha_{\De} \hat a_0\hat c_1^{m}$$
for each interval $\De\in {\mathcal S}_{\hat c_1}$ in $\Om$, and for each $x,y,z\in G_+$ with $x\not\perp u_m$ and $y\perp u_m$,
\medskip

(I) 
for uncountably many intervals $\Delta_{u}\in {\mathcal S}_{\hat u_m}$ with $\De_{u}\cap \De_{x'}\neq \emptyset$ for some $\De_{x'}\in {\mathcal S}_{\hat x}$, there is$\De_{y'}\in {\mathcal S}_{\hat y}$ 
such that $$\De_{u}(\hat s_m*\hat c_1^{-m}) \prec \De_{y'} \prec \De_{u},\; \hbox{and}$$

(II) 
for uncountably many intervals $\Delta_{u}\in {\mathcal S}_{\hat u_m}$ with $\De_{u}\cap \De_{x'}\neq \emptyset$ for some $\De_{x'}\in {\mathcal S}_{\hat x}$, there is  $\De_{y'}\in {\mathcal S}_{\hat y}$
such that $$\De_{u}(\hat s_m* \hat c_1^{-m}) \succ \De_{y'} \succ \De_{u}.$$

(III) if the supports of $\hat u_m$ and $\hat v_m$ are not disjoint, then we further have that  
$\Delta_u (\hat s_m*\hat c_1^{-m} )\;\cap \;\Delta_u\neq \emptyset$
for uncountably many $\Delta_u\in {\mathcal S}_{\hat u_m}$.
We can achieve this with $\hat s_m* \hat c_1^{-m}$ incomparable to the identity on any of these uncountably many $\De_u$ by the choice of the pertinent $\alpha_{\De}$ and $\{ h_{m,0}:m\in \N\}$.

\noindent We can also ensure that for all $m'\in \N$ and $f\in G_+$,
\medskip

(IV) for uncountably many $\Delta_{u}\in {\mathcal S}_{\hat u_m}$, there is $\De_{f_1}\in {\mathcal S}_{\hat f}$ such that $$\De_{u}(\hat s_m* \hat c_1^{-m}) \prec \De_{f_1*(\hat s_{m'}*\hat c_1^{-m'})} \prec \De_{u},\;\hbox{and}$$

(V)
for uncountably many $\Delta_{u}\in {\mathcal S}_{\hat u_m}$, there is $\De_{f_2}\in {\mathcal S}_{\hat f}$ such that $$\De_{u}(\hat s_m* \hat c_1^{-m}) \succ \De_{f_2*(\hat s_{m'}*\hat c_1^{-m'})} \succ \De_{u}.$$

But the $\ell$-subgroup of $G^{\dagger}$ generated by $G\cup  S^\prime$ is countable and $\Om$ is a minimal $\eta_1$-set. So we can further require that there is a subinterval $\De_0$ of $(\alpha_{\De}, \alpha_{\De}\hat a_0)$ disjoint from the convexification of $\bigcup_{n\in \N} supp(\hat g_n)$ such that  $\hat S^\prime$  maps $\De_{0}$ to itself and the restriction to $\De_{0}$ is a faithful representation of the free $\ell$-group on a countably infinite set $\{ z_n : n\in\N\}$ of generators under the map $z_n \mapsto \hat s_n * \hat c_1^{-n}$.
Hence we have an $\ell$-homomorphism ($x\mapsto \hat x$) of $G^{\dagger}$ into $B(\Om)$. 
As noted in Section 4, the restriction of this $\ell$-homomorphism to $G_{\flat}$ is  injective. The same is true for the restriction to $S$ by our construction.
Moreover, for any word $w$ in the alphabet of $G\cup \{a_0\}$ we have $[\hat w*\hat c_1^m,\hat s_{m'}]=1$ if $m,m'\in \N$ are distinct. 
Since $G^{\dagger}$ is countable, we can use the uncountability in (I) -- (V) to ensure that, in $\hat G^{\dagger}$, each $h_{0,m}$ acts as freely as possible on the $\ell$-subgroup of $B(\Om)$ generated by $\hat G\cup \{\hat a_0\}$.
That is, for any $\ell$-group word $w$ in the alphabet of $S\cup G\cup \{a_0\}$, if $w\not < 1$ in $G^{\dagger}$, we can find an interval on which 
$\hat w\not< 1$; similarly for $\not>$ and $\neq$.
So we can ensure that the $\ell$-homomorphism of $G^{\dagger}$ into $B(\Om)$ given by $x\mapsto \hat x$ is injective. Therefore, for such a choice of $\{ h_{m,0}:m\in \N\}$, we have

\begin{lemma}\label{iso}
With the above notation, 
$\hat G^{\dagger}\cong G^{\dagger}$. 
\end{lemma}

We are now ready to prove Proposition \ref{prop}.
\medskip

{\it Proof:}
It remains to show that $G^{\dagger}$ has soluble word problem.
\medskip

Fix a word $w$ in the alphabet of $G^{\dagger}$, say $w=\bigvee_I\bigwedge_J w_{i,j}$, where each $w_{i,j}$ is a group term.
So each $w_{i,j}$ is a group term in $\{ s_m:m\in \N\}\;\cup \;\{ a_0,c_1\}\;\cup \; \{g_n:n\in \N\}$.
Let $N_0$ be the $\ell$-subgroup of $G^{\dagger}$ generated by 
$S\;\cup\; G\;\cup \{ a_0\}$ and $N$ the $\ell$-ideal of $G^{\dagger}$ generated by $N_0$.
Now $G \ll a_0\ll c_1$, and $s_m \ll a_0 * c_1^m \ll c_1$ and $a_0 * c_1^{m+1} \perp a_0$ for all $m\in\N$.  By (\ref{orth}), $x\ll c_1$ for all $x\in N$. Hence 
$G^{\dagger}\cong N_0\wr (\langle c_1\rangle,\Z)$.
\medskip

Let $N_1$ be the $\ell$-subgroup of $G^{\dagger}$ generated by $G\cup S^\prime$.  Since $x\ll a_0$ for all $x\in G$ and $s_m\ll a_0 * c_1^m$ for all $m\in\N$, we get that $a_0\notin \bar N_1$ where $\bar N_1$ is the $\ell$-ideal of $G^{\dagger}$ generated by $G\cup S^\prime$. Then $G^{\dagger}/\bar N_1$ is $\ell$-isomorphic to the $\ell$-subgroup of $G^{\dagger}$ generated by $a_0, c_1$; this is $\ell$-isomorphic to $\langle a_0\rangle \wr (\langle c_1\rangle,\Z)$ and
$G^{\dagger}\cong N_1 \wr (\langle a_0\rangle, \Z) \wr (\langle c_1\rangle,\Z)$ . By Lemma \ref{Gwr}, it is enough to prove that $N_1$ has soluble word problem.
\medskip

Let $N_2$ be the $\ell$-ideal of $N_1$ generated by $G$ and $w$ be a word in the alphabet of $G\cup S^\prime$.
Say $w:=\bigvee_I \bigwedge_J w_{i,j}$, where each $w_{i,j}$ is a group word in this alphabet. Let $w^*$ be the result of replacing each occurrence of an element of $G$ by $1$.
\medskip

\noindent By the strengthing following (V), $N_1/N_2\cong S^\prime\cong F_{\aleph_0}$, the free $\ell$-group on a countably infinite set of generators.
Since $F_{\aleph_0}$ has soluble word problem, so does $N_1/N_2$;
and we can determine whether or not $w^*$  is $1$ in $N_1/N_2$.
 If $w^*\neq 1$ in this quotient, then its pre-image $w$ cannot be $1$ in $N_1$.
So assume that $w^*=1$; {\it i.e.}, $w\in N_2$.  
\medskip

Using the defining relations for $N_1$, we can effectively write each group word $w_{i,j}$ appearing in $w$ in the form $w_{i,j}^\prime w_{i,j}^{\prime\prime}$, where $w_{i,j}^\prime$ is a product of conjugates of elements of $G$ by elements of $S^\prime$, and $w_{i,j}^{\prime\prime}\in S^\prime$.  Moreover, $\bigvee_I\bigwedge_J \{w_{i,j}^{\prime\prime} : i\in I, j\in J\} = 1$ in $N_1$ since $w^*=1$.  
Now
\begin{equation}
\label{ij2}
w_{i,j}^\prime = \prod_k g_{i,j,k}*T_{i,j,k}(\{s_m*c_1^{-m}:m\in\N\}) 
\end{equation}

\noindent
with $g_{i,j,k}\in G$ and $T_{i,j,k}({\bf x})$ a group term involving a finite subset of variables from ${\bf x}$.
We may write
\begin{equation}
\label{T}
T_{i,j,k}(\{s_m*c_1^{-m}:m\in\N\})\;\;\hbox{as}\;\;(s_{m_{1,k}}*c_1^{-m_{1,k}})^{z_{1,k}}\dots (s_{m_{r,k}}*c_1^{-m_{r,k}})^{z_{r,k}},
\end{equation}
with $z_{1,k},\dots,z_{r,k}\in\Z\setminus\{0\}$ and $m_{1,k},\dots,m_{r,k}\in\N$ (not necessarily distinct).  

Assume first that $z_{1,k}>0$. Suppose that $g_{i,j,k}$ has subwords equal to $u_{m_{1,k}}$ or $u_{m_{1,k}}^{-1}$, say

\begin{equation}
\label{g}
g_{i,j,k}\;\;=\;\;g_{i,j,k,1}\;u_{m_{1,k}}^{\pm 1}\;g_{i,j,k,2}\;u_{m_{1,k}}^{\pm 1}\;\dots\;g_{i,j,k,\ell}
\end{equation}

\noindent
where each $g_{i,j,k,n}\in G$ contains no sub-occurrences of $u_{m_{1,k}}$ or $u_{m_{1,k}}^{-1}$ (to within equality in $G$) and may be 1 --- recall that we may (and have) assume(d) that all elements of $G$ are group words in the generators of $G$, and the solubility of the word problem for $G$ allows us to algorithmically determine the form (\ref{g}).  Replace the conjugate of $g_{i,j,k}$ appearing in $w_{i,j}^\prime$ by

\begin{equation}
\label{g2}
(g_{i,j,k,1}*T_{i,j,k})(v_{m_{1,k}}^{\pm 1}*T_{i,j,k}^{\prime})\dots 
(g_{i,j,k,\ell}*T_{i,j,k}),
\end{equation}
\noindent
where $$T_{i,j,k}^{\prime} = (s_{m_{1,k}}*c_1^{-m_{1,k}})^{z_{1,k}-1}(s_{m_{2,k}}*c_1^{-m_{2,k}})^{z_{2,k}}\dots (s_{m_{r,k}}*c_1^{-m_{r,k}})^{z_{r,k}}.$$

If $z_{1,k}>1$, determine if $v_{m_{1,k}}$ contains a subword equal in $G$ to $u_{m_{1,k}}^{\pm 1}$, and repeat the process with the subterms $v_{m_{1,k}}^{\pm 1}*T_{i,j,k}^{\prime}$.  Continue through at most $z_{1,k}$ steps to obtain a ``reduced" word with no further cancellation by just applying the relations (\ref{e4}) with $m=m_{1,k}$ and equality in $G$. In this way, we can reduce to the cases when $T_{i,j,k}^\prime$ begins with $(s_{m_{2,k}}*c_1^{-m_{2,k}})^{z_{2,k}}.$ 
\medskip

If $z_{1,k}<0$, write $g_{i,j,k}$ in the form (\ref{g2}) but with $v_{m_{1,k}}$ in place of $u_{m_{1,k}}$ and perform the same analysis interchanging $v_{m_{1,k}}$ and $u_{m_{1,k}}$. 
\medskip

We next consider if $v_{m_{1,k}}$ (or $u_{m_{1,k}}$ if $z_{1,k}<0$) contains a subword equal in $G$ to $u_{m_{2,k}}^{\pm 1}$ if $z_{2,k}>0$ (or $v_{m_{1,k}}^{\pm 1}$ if $z_{2,k}<0$) and repeat the process with $m_{2,k}$ in place of $m_{1,k}$. 
By continuing in this manner, we can write each $w_{i,j}^\prime$ in the form

\begin{equation} \label{last}
(g_{i,j,1}^\prime*t_{i,j,1})\dots
(g_{i,j,x}^\prime*t_{i,j,x})\cdot
(s_{m_1}^{\pm 1}c_1^{-m_1})\dots
(s_{m_k}^{\pm 1}c_1^{-m_k}),
\end{equation}
\noindent where $x,k\in\N,\;\;g_{i,j,1}^\prime,\dots,g_{i,j,x}^\prime\in G,\;\; m_1,\dots,m_k\in\N$, and all $t_{i,j,r}\in S^{\prime}$ are such that:

\noindent if $t_{i,j,r}$ begins $s_m*c_1^{-m}$, then $g_{i,j,r}^\prime$ contains no subword of form $u_m^{\pm 1}$, and 

\noindent if $t_{i,j,r}$ begins $s_m^{-1}*c_1^{-m}$, then $g_{i,j,r}^\prime$ contains no subword of the form $v_m^{\pm 1}$.
\medskip

The key to completing the algorithm to determine if the word $w$ ($\in N_2$) is not the identity in $N_1$ is the ``near-freeness" of the action of $\hat S^\prime$ on $(\alpha_\De, \alpha_\De \hat a_0)$.
We can now use the Holland-McCleary technique of Section 2.2 to determine if the representation of the resulting word is not the identity in $B(\Om)$. 
We need only consider the action on the interval $(\alpha_\De,\alpha_\De \hat a_0)$.
Extra considerations are needed, however. If $g\in G^+$, we must ensure that $\beta \hat g\geq \beta$ for all $\beta$ in all diagrams.
Similarly, if $f\leq g$ in $G$ and the action of $g$ already occurs in the diagram we can only extend the diagram so that $\beta f\leq \beta g$. 
And if $f\perp g$, then we must also have that $\beta f=\beta$ whenever $\beta g\neq \beta$.
For this reason, we must allow points to be fixed by elements of $\hat G\cup \hat S^\prime$ in our definition of legitimate diagrams. Also, if $1\leq g\leq u_m$, then $\beta\leq \beta(\hat g*(\hat s_m * \hat c_1^{-m}))\leq \beta\hat v_m$ for all $\beta\in (\alpha_{\De},\alpha_{\De}\hat a_0)$. We therefore take all the finitely many possibilities for $\la_0$ allowed by the word $w$: let $X_{i,j}$ be the set of all initial subwords of $w_{i,j}$ ($i\in I, j\in J$). For each subset $X_0$ of $X:=\bigcup \{ X_{i,j}:i\in I,\; j\in J\}$, let $X_0'$ denote its complement in $X$. Let $X_0\subseteq X$ be closed under initial subwords. Take any $\la_0\in \bigcup \{ supp(x): x\in X_0\} \setminus \bigcup \{ supp(y): x\in X_0'\}$ with all possible orderings (including equalities) for $\{ \la_0 u:u\;\hbox{an initial subword of}\; x\}$ (each $x\in X_0$). This provides a finite set of possibilities for $\la_0$ for each such subset of the finite set $X$. We proceed with each one that {\it is} consistent with the above considerations for $G\cup S^\prime$. We construct diagrams (as explained in Section 2.2) allowing all consistent possibilities for any $\hat s_{m_q}*\hat c_1^{-m_q}$ that appears in the resulting words as given in (\ref{last}) (according to (I) -- (V)), since $\hat h_m$ is locally arbitrary and can be positive, negative or a ``small" local perturbation ($m\in \Z_+$) for these $\lambda_0\in (\alpha_\De,\alpha_\De \hat a_0)$.  If, in any one of these finitely many consistent diagrams, we have $\lambda_0 \hat w \neq \lambda_0$, then $w\neq 1$ in $N_1$; and if $w\in N_2$, then such a legitimate diagram must exist if $w\neq 1$ in $N_1$.
[We illustrate with an example below.]

Thus we have an algorithm to determine if a word in the alphabet of $N_1$ is the identity or not. That is, $N_1$ has souble word problem, and hence so does  $G^{\dagger}$.
\pfend

{\bf Example.}
\medskip

\noindent Let $w=\bigvee_{i=1}^2\bigwedge_{j=1}^2 w_{i,j}\in N_1$, where $w_{1,1}=s_2*c_1^{-2},\;\;w_{1,2}=v_2^{-1},\;\;w_{2,1}=u_2(s_2*c_1^{-2})^{-1}$ and $w_{2,2}=g_3*(s_4*c_1^{-4})$.  Then $w^*=[(s_2*c_1^{-2})\wedge
1]\vee[(s_2*c_1^{-2})^{-1}\wedge 1]=1$, so $w\in N_2.$
Since $G$ has soluble word problem, we can determine whether or not $g_3\vee1=1$ in $G$ and whether or not $g_3=1$ in $G$.
\medskip

\noindent \underline{Case 1}$\;\;g_3\vee1\neq1$ in $G$.  
\medskip

\noindent By (IV), there are $\De_u\in {\mathcal S}_{\hat u_2}$ and $\De_g\in{\mathcal S}_{\hat g_3\vee 1}$ such that
$$\De_u(\hat s_2*\hat c_1^{-2})\prec \De_{g*(\hat s_4 *\hat c_1^{-4})} \prec \De_u.$$
Thus we obtain a legitimate diagram with $\lambda_0 \hat w > \lambda_0$ by taking $\lambda_0\in\De_{g*(\hat s_4 *\hat c_1^{-4})}$, since
$$\lambda_0<\lambda_0 \hat g_3*(\hat s_4*\hat c_1^{-4}) < \lambda_0\hat u_2(\hat s_2*\hat c_1^{-2})^{-1},\;{\mathrm so}$$
$$\lambda_0 \hat w = \lambda_0 \hat g_3*(\hat s_4 * \hat c_1^{-4}) > \lambda_0.$$

\noindent Hence there will be a consistent diagram with $\lambda_0<\lambda_0\hat w$, whence our algorithm will show $w\neq 1$ in $N_1$ if $g_3\vee1\neq1$ in $G$.
\bigskip

\noindent \underline{Case 2}$\;\;g_3\neq1=g_3\vee1$.  So $g_3<1$.
\medskip

\noindent By (IV), there are $\De_u\in{\mathcal S}_{\hat u_2}$ and $\De_g\in{\mathcal S}_{\hat g_3^{-1}\vee 1}$ such that
$$\De_u(\hat s_2* c_1^{-2})\prec \De_{g*(\hat s_4*\hat c_1^{-4})}\prec \De_u.$$
For $\lambda_0\in\De_{g*(\hat s_4 *\hat c_1^{-4})}$, we have
$$\lambda_0(\hat w_{2,1}\wedge\hat w_{2,2})\leq \lambda_0\hat w_{2,2}=\lambda_0(\hat g_3*(\hat s_4*\hat c_1^{-4}))<\lambda_0,\;{\mathrm and}$$ $$\lambda_0(\hat w_{1,1}\wedge\hat w_{1,2})\leq \lambda_0\hat w_{1,1}=\lambda_0(\hat s_2*\hat c_1^{-2})<\lambda_0.$$

\noindent Hence there will be a consistent diagram with $\lambda_0 \hat w<\lambda_0$.  Our algorithm will therefore display that $w\neq1$ in $N_1$ if $g_3<1$.
\bigskip

\noindent \underline{Case 3}$\;\;g_3=1$.
\medskip

\noindent Note that $\beta\hat w_{1,2}=\beta\hat v_2^{-1}\leq\beta$ and $\beta \hat w_{2,2}=\beta$ for all $\beta\in(\alpha_\De, \alpha_\De \hat a_0)$.
\medskip

\noindent Since $G$ has soluble word problem, we can determine whether or not $u_2\wedge v_2=1$ in $G$.
\medskip

\noindent \underline{Case 3(a)}$\;\;u_2\wedge v_2=1$ in $G$.
\medskip

\noindent By (II), there is $\De_u\in{\mathcal S}_{\hat u_2}$ such that
$$\De_u\prec\De_u(\hat s_2* \hat c_1^{-2}).$$
But $\De_u(\hat s_2* \hat c_1^{-2}) = \De_v$ for some $\De_v\in{\mathcal S}_{\hat v_2}$.
\medskip

\noindent For $\lambda_0\in\De_v$, we get
$\lambda_0\hat v_2^{-1}<\lambda_0$ and
$$\lambda_0\hat u_2(\hat s_2*\hat c_1^{-2})^{-1}=\lambda_0(\hat s_2*\hat c_1^{-2})^{-1}\in\De_u\prec\De_v.$$
Thus $\lambda_0\hat w = \lambda_0\hat v_2^{-1}<\lambda_0$.  So there is a consistent diagram with $\lambda_0\hat w<\lambda_0$, whence our algorithm will show that $w\neq1$ in $N_1$ in this case.
\medskip

\noindent \underline{Case 3(b)}$\;\;u_2\wedge v_2\neq1$ in $G$.
\medskip

\noindent By (II), there is $\De_u\in{\mathcal S}_{\hat u_2}$ and $\De_v\in{\mathcal S}_{\hat u_2\wedge\hat v_2}$ with
$$\De_u\prec\De_v\prec\De_u(\hat s_2*\hat c_1^{-2}).$$
For $\lambda_0\in\De_v$ we get $\lambda_0\hat v_2^{-1}<\lambda_0$ and
$$\lambda_0\hat u_2(\hat s_2*\hat c_1^{-2})^{-1}\in\De_v(\hat s_2*\hat c_1^{-2})^{-1}\prec\De_v.$$
Hence $\lambda_0\hat w=\lambda_0\hat v_2^{-1}<\lambda_0$ and again we have a legitimate diagram showing that $w\neq 1$ in $N_1$.
\bigskip

Therefore, in all circumstances:
$$g_3\vee1\neq1,\;\;g_3\vee1=1\neq g_3,\;\;{\mathrm and}\;\;(g_3=1 \;{\mathrm with\; either}\;u_2\perp v_2\;{\mathrm or}\;u_2\not\perp v_2)$$ our algorithm shows that $w\neq1$ in $N_1$.
\bigskip

If $w_0=w\vee1$, then our argument shows that $w_0\neq1$ in $N_1$ if $g_3\vee1\neq1$.  Since $w_{1,2},w_{2,2}\leq1$ if $g_3\vee1=1$, in all possible legitimate diagrams, $\lambda_0\hat w_0=\lambda_0$ if $g_3\vee1=1$.  So $w_0=1$ in $N_1$ if $g_3\vee1=1$.\pfend

\section{The proof of Theorem \ref{lbh}}

{\it Proof:} We can use Proposition \ref{prop} to define ${\frak S}(G)$ inductively.  Let $H$ be a recursively generated $\ell$-group with soluble word problem.  By Proposition \ref{prop}, there is a recursively generated $\ell$-group $H^{\dagger}$ with soluble word problem in which any two strictly positive elements of $H$ are conjugate.

Let $G(0):=G$ and $G(m+1):=G(m)^{\dagger}$.
Let ${\frak S}(G):= \bigcup \{ G(m):\; m\in \N\}$.
Then each $G(m)$ has a recursive set of generators by construction;
it has soluble word problem by Proposition \ref{prop} ($m\in \N$).
Thus the same is true of ${\frak S}(G)$ and any two strictly positive elements of ${\frak S}(G)$ are conjugate. Hence ${\frak S}(G)$ is $\ell$-simple, countable and has soluble word problem.
Therefore, ${\frak S}(G)$ has a recursive set of defining relations. By the proof of Theorem E of \cite{lH} (explained above in the second paragraph of Section 3) and Theorem \ref{lh}, ${\frak S}(G)$ can be $\ell$-embedded in a finitely presented $\ell$-group.
\pfend

\bigskip

{\bf Acknowledgement} 
\medskip

This paper is dedicated to Stephen H. McCleary on his $65^{th}$ birthday as a small thank you for many years of enjoyable collaboration.

\bigskip

Author's addresses:
\medskip

Queens' College,

Cambridge CB3 9ET,

England,
\medskip

and
\medskip

Department of Pure Mathematics and Mathematical Statistics,

Centre for Mathematical Sciences,

Wilberforce Rd.,

Cambridge CB3 0WB,

England
\medskip

amwg@dpmms.cam.ac.uk

\end{document}